\documentclass[12pt, a4paper]{article}
\usepackage{amsmath, amsthm,  amsfonts}

\def\ctr{{\rm ctr}}
\def\R{\mathbb R}
\def\d{\partial}
\def\P{{\cal P}}

\def\d{\partial}

\def\M{{\cal M}}
\def\V{{\cal V}}

\def\Sec{{\S}}
\def\Sp{{\cal S}}

\def\and{\hbox{\quad and \quad}}

\numberwithin{equation}{section}

\theoremstyle{plain}
\newtheorem{theorem}{Theorem}[section]
\newtheorem{lemma}{Lemma}[section]
\newtheorem{proposition}{Proposition}[section]

\theoremstyle{definition}

\theoremstyle{remark}

\title{The spherical phylon group and invariants of
the Laplace transform.}
\author{A.L. Carey$^1$, M.G. Eastwood$^1$, P.E. Jupp$^2$
and M.K. Murray$^1$}

\begin{document}
\maketitle

\begin{center}
$^1$Department of Pure Mathematics, University of Adelaide,
\linebreak
Adelaide SA 5005, Australia
\linebreak
$^2$School of Mathematical and Computational Sciences,
\linebreak
University of St Andrews, North Haugh, St Andrews KY16 9SS, U.K.
\end{center}
\medskip

\begin{abstract}
We introduce the spherical phylon group, a subgroup of
the
group of all formal diffeomorphisms of $\R^d$ that fix the
origin.
The   invariant theory of the spherical phylon group is used to
understand the invariants of the Laplace transform.
\end{abstract}

\section{Introduction}
\label{sec:intro}

There are various contexts in which it is of interest to find
 approximations for large $n$  to integrals of the form
\begin{equation}
\int_{\R^d}\exp\{-n f(x)\} b(x)  dx^1 \dots dx^d,
\label{eq:integral}
\end{equation}
where $f$ and $b$ are   smooth functions on $\R^d$ with $f$
having a Taylor series
at $0$ of the form
\begin{equation}
\frac{1}{2!} f_{ij}x^{i}x^{j} + \frac{1}{3!}
 f_{ijk}x^{i}x^{j}x^{k} +  \dots,
\label{eq:taylor_f}
\end{equation}
and $b(0) \neq 0$.
Here the matrix $(f_{ij})$ is required to be positive-definite and we assume
throughout the Einstein
summation convention.

Approximations to  integrals such as \eqref{eq:integral} are useful in
statistical
inference, see, for example  Barndorff-Nielsen \&  Wood (1995) and  Shun \&
McCullagh (1995).
Laplace's
method (see e.g. \Sec 6.4 of Bender \& Orszag (1978) or \Sec 8.3 of
Bleistein \&  Handelsman (1986)) provides an asymptotic expansion in $n$ of
\eqref{eq:integral} of the form
\begin{equation}
n^{-d/2}\{\Lambda_0(f, b) +  \Lambda_1(f, b) n^{-1/2} + \dots \} ,
\label{eq:asymp_exp}
\end{equation}
The coefficients $\Lambda_i$ depend on the function $f$ and the measure
$ bdx^1 \dots dx^d$.
The behaviour of \eqref{eq:asymp_exp} under reparameterisation has been
discussed in
Pace \& Salvan (1994) and Barndorff-Nielsen \&  Wood (1995).
Since the function $f$ and the measure
$b dx^1 \dots dx^d$ do not depend on any choice of coordinates on $\R^d$,
the integral \eqref{eq:integral} does not depend on any coordinate system.
It follows
that the coefficients $\Lambda _i (f,b)$ in \eqref{eq:asymp_exp} are
independent of the
choice of coordinates on $\R^d$.
In other words, the $\Lambda _i (f,b)$ are diffeomorphism invariants
of the pair $(f, b)$.
We show below that (if $d > 1$) these invariants are
essentially the only invariants of the pair $(f, b)$.
By this we mean that if we have two such pairs $(f, b)$ and $(g, c)$
such that $\Lambda _i (f,b) = \Lambda _i (g, c)$ for $i=0, 1, \dots$
then there is a diffeomorphism
mapping $(f, b)$ to  $(g, c)$.

There are a number of settings in which we could
attempt to prove this result. For example,
we could take $f$ and $b$ to be defined globally
 on a $d$-dimensional manifold
or locally in a
neighbourhood of the origin in $\R^d$. We consider here the case
that $f$ and $b$ are jets, or, because $\R^{d}$ has a canonical
co-ordinate system,  formal power series.
Since the asymptotic expansion \eqref{eq:asymp_exp} is obtained by Taylor
expansion
of $f$ and $b$ about 0, the coefficients $\Lambda _i(f,b)$ depend
on $f$ and $b$ only
through their Taylor expansions at $0$. Thus it is
possible to define the $\Lambda_i (f,b)$ when $f$ and $b$ are formal
power series. This approach is consistent with
Pace \& Salvan (1994) and fits in with previous work, on the
phylon group, of three of the authors contained in
 Barndorff-Nielsen {\it et al.} (1992), Carey \& Murray (1990)
 and Carey {\it et al.} (1996).
 Expressions
for the first few of the
$\Lambda_i$ are given in Pace \& Salvan (1994)
and, as might be expected, these
rapidly become complicated as $i$ increases. In the case that $b = 1$, a
combinatorial formula for the $\Lambda_i (f,b)$ can be obtained from
equation
(4) of  Shun \& McCullagh (1995).
We also present a general formula for the $\Lambda_i(f, b)$ in
\eqref{eq:formula}.
We find that the  $\Lambda_i(f, b)$ are polynomials in the
coefficients of $f$ and $b$ and rational in the square root of the
determinant of the Hessian matrix of $f$.
This is entirely analogous to classical invariant theory.

We begin with the calculation of the general formula for
the invariants and follow that with a proof of the
theorem.

\section{A general formula for the invariants}

 First let us fix some notation.
Given  $(f, b)$ let $\tilde f(x) = f(x) - (1/2)f_{ij}x^ix^j$ and
denote by $(f^{ij})$ the inverse of the Hessian matrix $(f_{ij})$.
Then we can write the integral \eqref{eq:integral} as
$$
\int_{\R^d} \exp\{-n \frac{1}{2}f_{ij}x^ix^j\}
\exp\{- n \tilde f(x)\} b(x)   dx^1 \dots dx^d.
$$
Standard techniques from distribution theory (see Appendix
\ref{sec:formula}) can be used to  show that if $i$ is odd
then $\Lambda_{i}(f, b) = 0$ and if $i$ is even then
\begin{equation}
\Lambda_i(f, b) =
\frac{(2\pi)^{d/2} }{\det(f_{ij})^{1/2}}
\sum_{ r+s-2l = i}  \frac{(-1)^l}{l!r!s!} (\tilde f^l)_{i_1 \dots
i_r}b_{j_{1}\dots
j_s} f^{k_1k_2}\dots f^{k_{r+s-1}k_{r+s}},
\label{eq:formula}
\end{equation}
where  $(1/r!)(\tilde f^l)_{i_1 \dots i_r}$
is the coefficient of $x^{i_{1}}\dots x^{i_{r}}$ in the Taylor
expansion of $(\tilde f)^{l}$, the $l$-th product of $\tilde f$,
and $k_1, k_2 |  \dots | k_{r+s - 1}, k_{r+s}$
runs through all partitions of
$i_1, \dots, i_s, j_1, \dots, j_s$ into subsets of size $2$.
For $i=0, 2$ the formula  \eqref{eq:formula} for the invariants
agrees with the results in Pace \& Salvan (1994).

\section{The  phylon group}
\label{sec:fps}
Denote by $\M(d)$ the set of all formal power series $f$,
in $d$ variables,
of the form \eqref{eq:taylor_f} with the matrix $(f_{ij})$
positive-definite and by $\V(d)$ the set of all
formal power series $b$ with $b(0) \neq 0$.
Denote by $\M^{k}(d)$ the set of all formal power
series $f$ which begin in the same way as \eqref{eq:taylor_f}
but which are finite of degree $k$. There is a
projection map
$$
j^{k}\colon \M(d) \to \M^{k}(d) ,
$$
the `$k$-jet' map which just truncates a
formal power series. Similarly we define $\V^{k}(d)$
and $j^{k}\colon \V(d) \to \V^{k}(d) $. Notice that $j^{0}(b)= b(0)$.

At various places we shall need one non-trivial fact
about formal power series. If $h$ is a smooth function
on a neighbourhood of $0$ in $\R^{d}$ then its Taylor
series  at $0$ defines a formal power series, which we denote by
$j^{\infty}(h)$.
 This map is well-known
not to be injective. However it follows  from a theorem
of Borel  (see e.g. p. 390 of Treves (1967)), that it is surjective,
that is every formal power series
is the Taylor series at $0$ of some smooth function. We shall
use this fact several times.

Denote by $\P(d)$ the set  of all formal
power series with values in $\R^d$ which have no
constant term and an invertible linear term. An element
$\psi$  in $\P(d)$ has the form
$\psi = (\psi^1, \dots, \psi^d)$ where
$$
\psi^{i}(x)  =\psi^i_j x^j + \frac{1}{2!} \psi^i_{jk} x^jx^k + \dots
$$
with $(\psi^i_j)$  an invertible matrix.
Under composition of  formal power series, $\P(d)$ forms a
group, called the {\em phylon group}. To prove this,
one can either construct inverses explicitly using the fact that
$(\psi^i_j)$ is an invertible matrix or replace $\psi$ by a smooth
function $\hat\psi$ using Borel's theorem and apply the inverse
function theorem.  (For further  details of $\P(d)$
and its representations  see
Barndorff-Nielsen {\it et al.} (1992),
Carey \& Murray (1990), Carey {\it et al.} (1996).)

We can also define  groups $\P^{k}(d)$ of finite formal
power series of degree $k$ with no constant term and with invertible
linear term. Their multiplication consists
of composition followed by truncation to order $k$.
With this multiplication the jet map
$j^{k}\colon \P(d) \to \P^{k}(d) $
is a homomorphism.
We denote the kernel of $j^{k}\colon \P(d) \to \P^{k}(d) $
 by $\P^{(k)}(d)$. This is  the
normal subgroup of all $\psi$ with components of the form
\begin{equation}
\psi^i(x) = \delta^i_j x^j + \frac{1}{ (k+1)!} \psi^i_{j_1\dots j_{k+1}}
 x^{j_1}\dots x^{j_{k+1}} + \dots .
\label{eq:condPk}
\end{equation}

Of particular interest is the first of these groups,
 $\P^{(1)}(d)$. The quotient $\P(d)/\P^{(1)}(d)$
is just $ GL(d)$, the general linear group of all
invertible $d$ by $d$ matrices.  Because a linear transformation
is a global diffeomorphism of $\R^d$, the general linear group
is also a subgroup of $\P(d)$ and it follows that
$\P(d)$ is a semi-direct product of $GL(d)$ and
$\P^{(1)}(d)$.

The phylon group acts on the space $\M(d)$ by $\psi$ sending $f$
to $\psi f$, where
\begin{equation}
(\psi f) (x) = f(\psi^{-1}(x)), \label{eq:f_action}
\end{equation}
and on $\V(d)$ by $\psi$ sending $b$ to $\psi b$, where
\begin{equation}
(\psi b) (x) = b(\psi^{-1}(x))
 \left| \det(\frac{\d (\psi^{-1})^i}{ \d x^j})\right|. \label{eq:b_action}
\end{equation}
Notice that a slight subtlety arises here in the interpretation
of the modulus signs. If
$$
a = a(0)  + a_i x^i + \dots
$$
is a power series with $a(0) \neq 0$ then $|a|$ denotes  the power
series defined by
\begin{equation}
|a|(x) = \begin{cases}
\phantom{-}a  & \text{if $a(0) > 0$ }\\
-a  & \text{if $ a(0) < 0$.}
\end{cases}
\label{eq:modulus}
\end{equation}

The group $\P^{k}(d)$ acts similarly on $\M^{k}(d)$ and
$\V^{k-1}(d)$.  These group actions
commute with the projections $j^{k}$ in the sense that
\begin{equation}
j^{k}(\psi f) = j^{k}(\psi)j^{k}(f) \quad\text{and}\quad
j^{k}(\psi b) = j^{k+1}(\psi)j^{k}(b). \label{eq:commute}
\end{equation}

We have noted in the
introduction that if we define the invariants $\Lambda_i$
 by the asymptotic expansion of \eqref{eq:integral}
 then they are manifestly invariant under
 local diffeomorphisms of $\R^{d}$ fixing the origin.
We need to know that when we consider the $\Lambda_{i}$
as functions on $\M(d) \times \V(d)$ they are
invariant under the action of $\P(d)$.
This can be proved using
Borel's theorem. We just  replace the formal $f$, $b$ and
$\psi$ by functions $\hat f$,  $\hat b$ and a local
diffeomorphism $\hat \psi$. Invariance
follows from the fact that $j^{\infty}(\hat \psi \hat f)
=\psi f$ and  $j^{\infty}(\hat \psi \hat b)
=\psi b$.

We wish to prove  the following theorem.
\begin{theorem}
Let $d > 1$. If $(f, b)$ and $(g, c)$ are
elements of $\M(d) \times \V(d)$ then
$$
\Lambda_i(f, b) = \Lambda_i(g, c) \quad \text{for all $i = 0, 1, 2,
\dots $}
$$
if and  only if there is a $\psi$ in $\P(d)$  such that
$$
(f, b) = (\psi g , \psi c).
$$
\label{th:phylon}
\end{theorem}

We begin our proof of the theorem by noting
that if $f$ is a function, not just a formal power series,
with Taylor series as in
\eqref{eq:taylor_f}  then Morse's Lemma (Milnor (1963)) shows that we can
change co-ordinates so that $f$ is the standard quadratic form
\begin{equation}
q(x) = \| x \|^2 \label{eq:quad}.
\end{equation}
By using Borel's theorem we can show that this
result is true also for formal power series. That means
that if $f \in \M(d)$ then we can find a $\psi$ in  $\P(d)$
such that $\psi f  = q$.

To understand when $(f, b)$ and $(g, c)$ can
be transformed one into the other by the phylon group,
we can
begin by transforming them so that $f$ and $g$ both become
$q$, where $q$ is defined in equation \eqref{eq:quad}.
Then we can try to transform $b$ into
$c$ without changing $q$. That is, we
act by the subgroup of
the phylon group that fixes  $q$.
We will call this the {\em spherical phylon
group} and denote it by $\Sp(d)$.  The spherical
phylon group is therefore the subgroup of all $\psi \in \P(d)$
such that
\begin{equation}
\|\psi(x)\|^2  = \|x\|^2.
\end{equation}

We can now see why the theorem is not going to work when
$d=1$.
In that case the spherical phylon group is just $O(1) = \{ +1 , -1\}$.
We discuss the one-dimensional case in Section \ref{sec:one-dim}.

Before discussing the proof of Theorem \ref{th:phylon}
we need to analyse the spherical phylon group in more detail.

Notice that the intersection of $GL(d)$ with $\Sp(d)$
is precisely the group $O(d)$ of orthogonal transformations.
Moreover, if $\psi$, defined by
$$
\psi(x) = \psi^i_j x^ j + \frac{1}{ 2!} \psi^i_{jk} x^j x^k + \dots ,
$$
is in the phylon group then it is easy to check that $X =  \psi^i_j$
is orthogonal. Hence $X^{-1}\psi$ is also in the spherical
phylon group, so that $\Sp(d)$ is the semi-direct product
of $O(d)$ and $\Sp^{(1)}(d)$, where
$\Sp^{(1)}(d) = \Sp(d) \cap \P^{(1)}(d)$.
If we apply the condition for an element $\psi$ with components
$$
\psi^i(x) = \delta^i_j x^ j + \frac{1}{ 2!} \psi^i_{jk} x^j x^k + \dots
$$
to be in $\Sp(d)$, we see that this is so if
\begin{equation}
2\delta_{i(j_1}\psi^i_{j_2\dots j_{k+2})} + \delta_{ij}
\sum_{r=2}^{k+1} \psi^i_{(j_1\dots j_r}\psi^j_{j_{r+1}\dots j_{k+2)} }= 0
\label{eq:sph-cond}
\end{equation}
for all $k = 1, 2, \dots $, where we use the
convention that indices appearing between
round brackets are symmetrised over.
Define $\Sp^{(k)}(d) = \Sp(d) \cap \P^{(k)}(d)$.
If $\psi$ is an element of  $\Sp^{(k)}(d)$ then it satisfies
 \eqref{eq:condPk} and so the condition
in equation \eqref{eq:sph-cond} reduces to
\begin{equation}
\delta_{i(j_1}\psi^i_{j_2\dots j_{k+2})} = 0.
\label{eq:sph-condk}
\end{equation}

Let us write $\Lambda_i(b)$ for $\Lambda_i(q, b)$.
Again,
straightforward calculations using moments of the
multivariate normal distributions or direct application
of \eqref{eq:formula} shows that
\begin{equation}
\Lambda_i(b) = \pi^{d/2} \frac {1} {(i/2)!2^{i}} \ctr_{i}(b),
\label{eq:red-inv}
\end{equation}
where the {\em i-th complete trace}, $\ctr_{i}$, of $b$ is defined by
\begin{equation}
\ctr_{i}(b)=
\left\{
\begin{array}{ll}
0 & \mbox{if $i$ is odd} \\
\sum_{k_{1}\dots k_{j} } b_{k_{1}k_{1}\dots k_{j}k_{j}}
& \mbox{if $i=2j$ is even}
\end{array}
\right.
\label{eq:ctr}
\end{equation}
and $(1/r!)b_{{i_1}\dots{i_r}}$ is the
coefficient of $x^{i_1}\dots x^{i_r}$ in the
Taylor expansion of $b$.

Theorem \ref{th:phylon} is  a consequence of the following
theorem:
\begin{theorem}
If $b, c$ are elements of $\V(d)$  with
$d > 1$ then there is an element
of the spherical phylon group $\psi$ such that $\psi(b) = c$
if and only if $\Lambda_i(b) = \Lambda_i(c)$ for all $i$.
\label{th:spherical}
\end{theorem}

\section {Proof of  theorem \ref{th:spherical}}
In one direction the proof follows from the invariance
of the $\Lambda_i$  under the action of the phylon group.
To prove the other direction, start with $b$ and $c$ in
$\V(d)$ such that $\Lambda_{i}(c) = \Lambda_{i}(b)$
for all $i = 0, 1, \dots$ .
We shall prove below the following Proposition:

\begin{proposition}
If $b$ and $c$ are elements of $ \V(d)$ such that
$\Lambda_{i}(c) = \Lambda_{i}(b)$ for all $i$ and
$j^{k-1}(c) = j^{k-1}(b)$ then there exists $\psi_{k} \in \Sp^{(k)}(d)$
such that $j^{k}(c) = j^{k}(\psi_{k}b)$.
\label{prop}
\end{proposition}

Before proving Proposition \ref{prop}, let us see how it proves
Theorem \ref{th:spherical}. It follows from \eqref{eq:red-inv}
that $c(0) = b(0)$ so, by definition,  $j^0(c) = j^0(b)$.  We now use
Proposition \ref{prop}
to manufacture a sequence  of elements  $\psi_{k} \in \Sp^{(k)}(d)$ for
each $k = 1, 2, \dots$  such that
$$
j^{k}(c) = j^{k}(\psi_{k}\psi_{k-1} \dots \psi_{1}b).
$$
Note that when we multiply $\psi_{k-1} \dots \psi_{1}$
by $\psi_{k}$ on the left we have
\begin{align*}
j^{k}(\psi_{k}\psi_{k-1} \dots \psi_{1}) &=
j^{k}(\psi_{k})j^{k}(\psi_{k-1} \dots \psi_{1})\\
&=j^{k}(\psi_{k-1} \dots \psi_{1}),\\
\end{align*}
because $\psi_{k} \in \Sp^{(k)}(d)  =
\text{ker}(j^{k})$. Hence multiplying
$\psi_{k-1} \dots \psi_{1}$ by $\psi_{k}$ changes
no terms of degree lower than $k+1$. This means
that there is a well-defined formal power series
$\psi \in \P(d)$ such that for all $k= 1, 2, \dots$
$$
j^{k}(\psi) = j^{k}(\psi_{k-1} \dots \psi_{1}).
$$
Then we have
\begin{align*}
j^{k}(c) &= j^{k}(\psi_{k}\psi_{k-1} \dots \psi_{1}b)\\
&=j^{k+1}(\psi_{k}\psi_{k-1} \dots \psi_{1})j^{k}(b)\\
&= j^{k+1}(\psi)j^{k}(b)\\
&=j^{k}(\psi b)\\
\end{align*}
for all $k = 1, 2, \dots$.  Hence $c = \psi b$, proving
Theorem \ref{th:spherical}.

It remains to prove Proposition \ref{prop}.
Write
$$
b = \alpha + B + \tilde b
$$
and
$$
c = \alpha + C + \tilde c,
$$
where $\alpha$ is a polynomial of degree of $k-1$,  $B$
and $C$ are homogeneous of degree $k$ and $j^k(\tilde b) = 0 =
j^k(\tilde c)$.

For any $\psi$ in $ \Sp^{(k)}(d)$, write
$$
\psi^i(x) = \delta^i_j x^j +
\frac{1}{ (k+1)!}\psi^i_{j_1\dots j_{k+1}} x^{j_1} \dots x^{j_{k+1}}
     + O(|x|^{k+2}).
$$
Then
$$
\frac{\partial \psi^i }{\partial x^j} =
\delta^i_j  + \frac{1}{k!}\psi^i_{jj_1\dots j_{k}} x^{j_1} \dots x^{j_{k}} +
O(|x|^{k+1}),
$$
so that
\begin{equation}
\det \left(\frac{\partial \psi^i }{\partial x^j}\right) =
1 + \frac{1}{k!}\psi^i_{ij_1\dots j_{k}} x^{j_1} \dots x^{j_{k}} +
O(|x|^{k+1}).
\label{eq:jac}
\end{equation}
Note that because the power series in
\eqref{eq:jac} is $1$ at $0$, its absolute
value is itself (by \eqref{eq:modulus}). Then  from \eqref{eq:b_action} we
have
$$
\psi^{-1}b = (b\circ \psi )\det
\left(\frac{\partial \psi^i }{\partial x^j}\right),
$$
so that
$$
(\psi^{-1}b)(x) = \alpha(x) + B(x) + \frac{b(0)}{k!}\psi^i_{ij_1\dots
j_{k}} x^{j_1} \dots
x^{j_{k}} + O(|x|^{k+1}).
$$
 If  we can find a $\psi^i_{j_1\dots j_{k+1}}$ satisfying
 \eqref{eq:sph-condk}  such that
\begin{equation}
C(x)  = B(x) + b(0)\psi^i_{ij_1\dots j_{k}}
x^{j_1} \dots x^{j_{k}}  \label{eq:tensor}
\end{equation}
then defining $\psi$ by
$$
\psi^i(x) = \delta^i_j x^j +
\frac{1}{ (k+1)!}\psi^i_{j_1\dots j_{k+1}} x^{j_1} \dots x^{j_{k+1}}
$$
gives
$$
j^{k}(b)= j^{k}(\psi c),
$$
which proves Proposition \ref{prop}.

Thus  we have reduced the proof of Proposition \ref{prop} (hence that
of  Theorem \ref{th:spherical}) to solving  \eqref{eq:tensor}.
Denote by $K$ the subspace
of $\R^d \otimes S^{k+1}(\R^d)$ of tensors $b_{j_1\dots j_{k+2}}$
such that
$$
b_{(j_1j_2\dots j_{k+2})} = 0.
$$
Here, and below, $S^{k+1}(\R^d)$ denotes  the space of totally symmetric
tensors of order $k+1$. Define a map  $T \colon K \to S^{k}(\R^d)$
 by
$$
T(b_{ij{j_1}\dots {j_{k}}})  = b_{ii j_1 \dots j_{k} }.
$$
Recall the definition of the complete trace  from equation
\eqref{eq:ctr}. To show that $\eqref{eq:tensor}$ can be solved
for $\psi^i_{j_1\dots j_{k+1}}$, we use \eqref{eq:red-inv} and the
following Lemma.
\begin{lemma}
An element $b \in S^{k}(\R^d)$ is in the image
of $T$ if and only if $\ctr_{k}(b) = 0$.
\label{lem:one}
\end{lemma}
\begin{proof}
Calculation shows that $\ctr_{k} \circ T = 0$. To show that the kernel of
$\ctr_{k}$ is contained in the image of $T$, we shall use a decomposition
\eqref{eq:dec} of $S^k(\R^d)$.
For any symmetric product $S^q(\R^d)$ we define $S^{[q]}(\R^d)$ as the set
of trace-free tensors, that is,
$$
S^{[q]}(\R^d) = \{ Z_{{j_1}\dots j_q} \in S^q(\R^d) \vert
\delta^{ij} Z_{ij{j_1}\dots j_{q-2}} = 0 \} .
$$
If $k-q = 2m$ then we can find a copy of $S^{[q]}(\R^d)$ in
$S^{k}(\R^d)$ as the image of the map
$$
W_{{i_1}\dots{i_q}} \mapsto \delta_{(j_1 j_2}\dots
                   \delta_{j_{2m-1}j_{2m}}W_{{j_{2m+1}}\dots{j_{k}})}.
$$
A calculation (or the corresponding decomposition of $O(d)$-modules
given in \Sec 19.5 and Lemma 17.15 of
Fulton \& Harris (1991)) shows that $S^{k}(\R^d)$
is the direct sum of all such images.
Thus we have
\begin{equation}
S^{k}(\R^d) = S^{[k]}(\R^d) \oplus  S^{[k-2]}(\R^d) \oplus \dots ,
\label{eq:dec}
\end{equation}
ending in $S^{[1]}(\R^d)$ or $S^{[0]}(\R^d)$, depending on
whether $k$ is odd or even. If $b \in S^{k}(\R^d)$ then
the complete trace of $b$ is either zero (if $k$ is odd) or
a non-zero multiple of the component of $b$ in $S^{[0]}(\R^d)$
(if $k$ is even).

It suffices to show that every tensor in
the irreducible component $S^{[q]}(\R^d)$ for
$0 < q \leq k$ is in the image of the map $T$.
Let
$$
Z_{j_1\dots j_{k}} = \delta_{(j_1 j_2}\dots
                   \delta_{j_{2m-1}j_{2m}}W_{{j_{2m+1}}\dots{j_{k}})}
$$
be  such a tensor.  Then consider
$$
X_{ii_1\dots i_{k+1}} =
\delta_{i(i_1}\dots  \delta_{i_{2m}i_{2m+1}}Z_{i_{2m+2}\dots i_{k+1})}
-\delta_{(i_1 i_2}\dots  \delta_{i_{2m+1}i_{2m+2}}Z_{i_{2m+3}\dots
i_{k+1})i}.
$$
This clearly satisfies
$$
X_{(ii_1\dots i_{k+1})}= 0.
$$
A calculation shows that
$$
T( X_{ii_1 \dots i_{k+1}}) =
\frac{d+ q -2 }{ k+1} Z_{i_1 \dots i_{k}}.
$$
Since $d > 1$, this proves the result.
\end{proof}

\section{The one-dimensional case}
\label{sec:one-dim}

We have already noted that Theorem \ref{th:phylon} cannot be
true in the case that $d=1$. To be able to
determine when there is a $\psi$ such that
$\psi(f, b) = (g, c)$ we must either restrict
the possible $(f, b)$ or define  extra invariants.
We consider both possibilities in this section.

The phylon group $\P(d)$ has two connected components.
If $\psi$ is an element of $\P(d)$ with components
$$
\psi^{i}(x)  =\psi^i_j x^j + \frac{1}{2!} \psi^i_{jk} x^jx^k + \dots
$$
then the different connected components of $\P(d)$ are determined
by the sign of the determinant of the matrix $(\psi^i_j)$.
Denote by $\P_0(d)$ the connected component of $\P(d)$
containing the identity. This is the component where
$\det(\psi^i_j) > 0$.

In the case that $d=1$ we shall write a general element
$\psi \in \P(1)$  using  different notation as
$$
\psi(x)  =\psi_1 x + \frac{1}{2!} \psi_2 x^2 + \dots .
$$
Then $\psi$ is in $\P_0(1)$ if
$\psi_1 > 0$. The important fact about this
group for our purposes is
\begin{theorem} If $f$ and $g$ are elements of $\M(1)$ then there
is a unique $\psi \in \P_0(1)$ such that $\psi f = g $.
\label{th:unique}
\end{theorem}
\begin{proof} Morse's Lemma shows that $\psi$ exists. To show that
$\psi$ is unique  it is enough to show that if $\phi q = q$ for
some $\phi \in \P_0(1)$
then $\phi = 1$. If $\phi q = q$ then
$\phi$ preserves length and, as we have noted previously
in one dimension, this means that $\phi(x) = \pm x$ for all $x$.
For $\phi$  to be in the connected component of the identity we must have
$\phi(x) = x$, so that $\phi = 1$.
\end{proof}

Notice that in the one-dimensional case it follows
from  \eqref{eq:red-inv} and \eqref{eq:ctr} that
for any $b \in V(1) $ the invariant $\Lambda_i(b) $ is
a non-zero multiple of
the $i$th coefficient in the Taylor expansion
$$
b(x) = b_0 + b_1 x + \frac{1}{2!} b_2 x^2 + \dots ,
$$
if $i$ is even and is zero if $i$ is odd.

The definitions of even and odd functions extend also to
power series. We define a power series  $f$ to be
 {\em even} if $f(-x) = f(x)$
for all $x$ and {\em odd} if $f(-x) = - f(x)$ for all $x$. If $g$ is
any power series we define
$$
g^+(x) = \frac{1}{2}(g(x) + g(-x)) \quad\text{and}\quad
g^-(x) = \frac{1}{2}(g(x) - g(-x)).
$$
Then $g(x) = g^+(x) + g^-(x)$ is the unique expansion
of $g(x)$ as a sum of an even power series  and an odd power series.
An even (odd)  power series has only even (odd) terms in its Taylor
expansion. Consider the  subgroup $P_{\text{odd}}(1)$ of $P_0(1)$
consisting of all $\psi$ which are odd.  Then  we have
\begin{theorem}
Suppose that $(f, b)$ and $(g, c)$ are elements of  $\M(1)\times  \V(1)$.
with $f$ and $g$ even  and that
$\Lambda_i(f, b) = \Lambda_i(g, c)$ for all $i=0, 1, 2, \dots $.
Then there is a $\psi \in \P_{\text{odd}}(1)$ such that
$\psi f = g$ and $\psi b^+ = c^+$.
\label{th:eq-case}
\end{theorem}
\begin{proof}
{}From Theorem \ref{th:unique} we have that there is a unique $\psi_1 \in
\P_0(1)$ satisfying $\psi_1 f = q$.
 Define $\tilde \psi_1 (x) = -\psi_1(-x)$.
Then $\tilde \psi _1 f(x) = f(\tilde \psi_1^{-1} (x)) =
f(-\psi_1^{-1} (-x) ) = f(\psi_1^{-1}(-x)) = q(-x)  = q(x)$. Hence
by uniqueness $\tilde \psi_1= \psi_1$,
 so that $\psi _1 \in \P_{\text{odd}}(1)$.
Similarly there is a $\psi_2 \in \P_{\text{odd}}(1)$ such that
$\psi_2 g = q$.

By invariance we deduce that $\Lambda_i(\psi_1 b) = \Lambda_i(\psi_2
c)$ and hence
that $\psi_1 b$ and $\psi_2 c$ have the same even Taylor
coefficients.  Thus $(\psi_1 b)^+ = (\psi_2 c)^+$. Consider
now a   $b$ of the form $b(x) = x^k$. Letting $d\psi_1$
denote the derivative of $\psi_1$ we have
$$
(\psi_1 b)(x)  = d\psi_1^{-1}(x) (\psi_1^{-1}(x))^k ,
$$
so that
$$
(\psi_1 b)(-x) = (-1)^k (\psi_1 b)(x).
$$
It follows that for any $b$, because $b^+$ is even,
we must have that $\psi_1 b^+ $ is even and similarly
$\psi_1 b^-$ is odd, so that
$$
\psi_1 b = \psi_1 b^+ + \psi_1 b^-
$$
is the decomposition of $\psi_1 b$ into even and
odd functions and, in particular,
 $(\psi_1 b)^+ = \psi_1 b^+ $.  Similarly
$(\psi_2 c)^+ = \psi_2 c^+ $ so we have
$ \psi_1 b^+  =  \psi_2 c^+ $. Letting $\psi = \psi_2^{-1} \psi_1
$ gives the desired result.
\end{proof}

Consider now the case of general   $(f, b)$ and $(g, c)$
with
$\Lambda_i(f, b) = \Lambda_i(g, c)$ for all $i = 0, 1, 2, \dots$
We proceed as  in the proof of Theorem \ref{th:eq-case}.
We  can find  $\psi_1$ and
$\psi_2$ such that $\psi_1 f = q$ and $\psi_2 g = q$.  We then
have
$\Lambda_i(\psi_1 b) = \Lambda_i( \psi_2 c) $ for all $i =0,  1, 2,
\dots $ and hence the even Taylor coefficients of
$\psi_1 b $ and $\psi_2 c $ equal but we have no information about
the odd coefficients.  If $(f, b) $  is in $\M(1) \times \V(1)$
we define $\lambda_i(f, b) $ to be the $i$th Taylor coefficient
of $\psi b$,  where $\psi \in \P_{\text{odd}}(1)$ and
$\psi f = q$. Then we have
\begin{theorem}
If $(f, b)$ and $(g, c)$ are
elements of $\M(1) \times \V(1)$ then
$$
\lambda_i(f, b) = \lambda_i(g, c) \quad \text{for all $i = 0, 1, 2,
\dots $}
$$
if and  only if there is a $\psi$ in $\P_0(1)$  such that
$$
(f, b) = \psi(g , c).
$$
\label{th:one-phylon}
\end{theorem}
\begin{proof}
The proof follows from the discussion we have already
given that motivated the definition of the $\lambda_i$.
\end{proof}

We have no explicit formula for the
calculation of the   $\lambda_i(f, b)$. Let us
show that an algorithm is possible. This is clear,
except for finding the $\psi$ that solves $\psi f = q$. This
equation is equivalent to
$f(\psi^{-1}(x)) = x^2$  and hence  $f(x) = \psi(x)^2 $.
If we let
$$
f (x) = \frac{1}{2!}f_2  x^2 + \frac{1}{3!} f_3 x^3 + \dots
$$
and
$$
\psi(x) = \psi_1 x + \frac{1}{2!} \psi_2 x^2 + \dots
$$
then the equation $f(x) = \psi(x)^2$ is equivalent to
an infinite system of equations
\begin{equation}
\begin{split}
\frac{1}{2!} f_2 &=  \psi_1^2 \\
\frac{1}{3!} f_3 &=  \psi_1\psi_2 \\
\frac{1}{4!} f_4 &= \frac{1}{4} \psi^2_2 + \frac{1}{3} \psi_1  \psi_3 , \\
& \vdots
\end{split}
\end{equation}
which we can solve recursively for $\psi_1, \psi_2, \dots $.

\begin{appendix}

\renewcommand{\theequation}{A.\arabic{equation}}
\makeatletter
\@addtoreset{equation}{section}
\setcounter{equation}{0}
\makeatother

\section*{Appendix A: Moment generating functions}
\label{sec:formula}
We include here some standard results on moment generating functions
which are
needed to derive formula \eqref{eq:formula}. The normal distribution
with mean  $0$ and variance matrix $\Sigma = (\Sigma^{ij})$ on $ \R^{d}$
has probability density function
\begin{equation}
\phi(x) = \frac{1}{(2\pi)^{d/2}(\det\Sigma)^{1/2}}
\exp\{-\frac{1}{2} \Sigma_{ij} x^{i} x^{j}\} ,
\label{eq:normal}
\end{equation}
where we assume the summation convention
 and the matrix $(\Sigma_{ij})$ is the
inverse of the matrix $ (\Sigma^{ij})$.  If $f$ is a function,
we define its expectation with respect to this distribution by
$$
E[f(x)]  = \int_{\R^d} f(x) \phi(x) dx^{1}\dots dx^{d}.
$$

The moment generating function $M(t)$, where $t = (t_1, \dots, t_d)$,
is defined by
\begin{equation}
M(t) = E[\exp(t_{i}x^{i})].
\label{eq:moment}
\end{equation}
By completing the square in \eqref{eq:moment} and using the fact that
the integral of \eqref{eq:normal} over all of of $\R^{d}$ is $1$, we
can show that
\begin{equation}
M(t) = \exp(\frac{1}{2} \Sigma^{ij} t_{i} t_{j}).
\label{eq:moment-normal}
\end{equation}

Repeated differentiation of  \eqref{eq:moment-normal}
followed by evaluation at $t=0$
shows that
$$
E(x^{i_1}x^{i_2}\dots x^{i_k}) = \sum \Sigma^{k_1 k_2} \Sigma^{k_3
k_4} \dots \Sigma^{k_{r-1} k_r}
$$
where the sum is over all partitions of $i_1, \dots, i_r$
into subsets of size $2$.
Hence if $r$ is odd  this expectation is zero.
It is straightforward now to calculate the formula in
\eqref{eq:formula}.

\end{appendix}
\vskip 12pt

\noindent{\bf Acknowledgement:} The support of the Australian
Research Council is  acknowledged.  PEJ thanks the other authors
for their hospitality during his visit to  Adelaide.

\section*{References}

\begin{description}
\item[]
Barndorff-Nielsen, O.E., Bl\ae sild, P.,  Carey, A.L.,
Jupp, P.E.,   Mora, M. \& M.K. Murray, M.K. 1992
Finite dimensional algebraic representations of the
infinite phylon group.
{\em Acta. Appl. Math.} {\bf 28}, 219--252.

\item[]
Barndorff-Nielsen, O.E. \& Wood, A.T.A. 1995
On large deviations and choice of ancillary for $p^*$ and the modified
directed likelihood.
 Research Report {\bf 299}, Dept. of Theoretical
Statistics, Aarhus University.

\item[]
Bender, C.M.  \& Orszag, S.A. 1978
{\em Advanced mathematical methods for scientists and engineers}.
Singapore: McGraw-Hill.

\item[]
Bleistein, N. and Handelsman, R. 1986
{\em Asymptotic expansions of integrals.}
 New York:       Dover.

 \item[]
Carey, A.L.  \&  Murray, M.K. 1990
Higher-order tensors, strings and new tensors.
{\em Proc. R. Soc. Lond. A.}
{\bf 430}, 423-432.

\item[]
Carey, A.L.,  Jupp, P.E.  \& Murray, M.K. 1996
The phylon group and statistics.
In  {\em Algebraic Groups and Related Subjects; a
Volume in Honour of R. W. Richardson\/}
Australian Mathematical Society Lecture Series,
(ed. G. I. Lehrer {\it et al.}),
pp.~46--60.
Cambridge: Cambridge University Press.

\item[]
Fulton, W.  \&  Harris, J. 1991
{\em Representation Theory, A First Course.}
New York: Springer-Verlag.

\item[]
Milnor, J. 1963
{\em Morse Theory.}
Princeton: Princeton University Press.

\item[]
Pace, L.  and  Salvan, A. 1994
A note on invariance of the Laplace approximation  under change of variable.
 Research Report 1994.7,
Department  of Statistical Sciences, University of Padova.

\item[]
Shun, Z.  \&  McCullagh, P. 1995
Laplace approximation of high dimensional integrals.
{\em J. Roy. Statist. Soc. B}
{\bf 57}, pp.~749--760.

\item[]
Treves, F. 1967
{\em Topological vector spaces, distributions and kernels.}
New York: Academic Press.
\end{description}

\end{document}